\documentclass[10pt,twoside]{article}
\usepackage{graphicx}
\usepackage{amsmath}
\usepackage{Latex-document}

\title{\bf  ``Algebraic Truths'' \vskip -2mm \textit{vs} \vskip -2mm ``Geometric Fantasies'':
\vskip -2mm Weierstrass' Response to Riemann\vskip 6mm}

\author{U. Bottazzini\vspace*{-0.5cm}\thanks{Dipartimento di matematica,
Universit\`a di Palermo, via Archirafi 34, 90123 Palermo, Italy.
E-mail: bottazzi@math.unipa.it}}
\date{\vspace{-8mm}}

\begin{document}
\maketitle

\thispagestyle{first} \setcounter{page}{923}

\begin{abstract}\vskip 3mm

In the 1850s Weierstrass succeeded in solving the Jacobi
inversion\ problem for the hyper-elliptic case, and claimed he was
able to solve the general problem. At about the same time Riemann
successfully applied the geometric methods that he set up in his
thesis (1851) to the study of Abelian integrals, and the solution
of Jacobi inversion problem. In response to Riemann's
achievements, by the early 1860s Weierstrass began to build the
theory of analytic functions in a systematic way on arithmetical
foundations, and to present it in his lectures. According to
Weierstrass, this theory provided the foundations of the whole of
both elliptic and Abelian function theory, the latter being the
ultimate goal of his mathematical work. Riemann's theory of
complex functions seems to have been the background of
Weierstrass's work and lectures. Weierstrass' unpublished
correspondence with his former student Schwarz provides strong
evidence of this. Many of Weierstrass' results, including his
example of a continuous non-differentiable function as well as his
counter-example to Dirichlet principle, were motivated by his
criticism of Riemann's methods, and his distrust in Riemann's
``geometric fantasies''. Instead, he chose the power series
approach because of his conviction that the theory of analytic
functions had to be founded on simple ``algebraic truths". Even
though Weierstrass failed to build a satisfactory theory of
functions of several complex variables, the contradiction between
his and Riemann's geometric approach remained effective until the
early decades of  the 20$^{th}$ century.

\vskip 4.5mm


\end{abstract}

\vskip 12mm

\section*{Introduction} \label{introduction}\setzero
\vskip-5mm \hspace{5mm }

In 1854 Crelle's {\it Journal} published a paper on Abelian
functions by an unknown school teacher. This paper announced the
entry in the mathematical world of a major figure, Karl
Weierstrass (1815-1897), who was to dominate the scene for the
next forty years to come. His paper presented a solution of Jacobi
inversion problem in the hyper-elliptic case. In analogy with the
inversion of elliptic integrals of the first kind, Jacobi
unsuccessfully attempted a direct inversion of a hyper-elliptic
integral of the first kind. This led him to consider multi-valued,
``unreasonable'' functions having a ``strong multiplicity'' of
periods, including periods of arbitrarily small (non-zero)
absolute value. Jacobi confessed he was ``almost in despair''
about the possibility of the inversion when he realized ``by
divination'' that Abel's theorem provided him with the key for
resurrecting the analogy with the inversion of elliptic integrals
by considering the sum of a suitable number of (linearly
independent) hyper-elliptic integrals instead of a single
integral. In his memoir submitted to the Paris Academy in 1826
(and published only in 1841) Abel had stated a theorem which
extended Euler's addition theorem for elliptic integrals to more
general (Abelian) integrals of the form $\int R(x,y) dx$  in which
$R(x,y)$ is a rational function and $y= y(x)$ is an algebraic
function defined by a (irreducible) polynomial equation $f(x,y) =
0$. According to Abel's theorem, the sum of any number of such
integrals reduces to the sum of a number $p$ of linearly
independent integrals and of an algebraic-logarithmic expression
({\it p} was later called by Clebsch the genus of the algebraic
curve $f(x,y) = 0$). In 1828 Abel published an excerpt of his
Paris memoir dealing with the particular (hyper-elliptic) case of
the theorem, when $f(x,y)=y^2 -P(x)$, $P$ is a polynomial of
degree $n > 4$ having no multiple roots. In this case $p =
[(n-1)/2]$, and for hyper-elliptic integrals of the first kind
$\int \frac{Q(x)dx}{\sqrt{P(x)}}$ ($Q$ is a polynomial of degree
$\leq p-1$) the algebraic-logarithmic expression vanishes
(\cite{Abel1881}, vol. 1, 444-456).

On the basis of Abel's theorem in 1832 Jacobi formulated the
problem of investigating the inversion of a system of $p$
hyper-elliptic integrals

$$u_k=\sum_{j=0}^{p-1}\int_0^{x_j}\frac{x^k dx}{\sqrt{P(x)}}     \mbox{     }   \mbox{     }  (0\leq k \leq p-1)
  \mbox{     }  \mbox{     }   (\textrm{deg}P = 2p+1 \mbox{ or}  \mbox{     }   2p+2) $$
by studying $x_0,x_1,\cdots,x_{p-1}$ as functions of the variables $u_0,u_1,\cdots,u_{p-1}$. These functions
$x_i=\lambda_i(u_0,u_1,\cdots,u_{p-1})$ generalized the elliptic functions to $2p$-periodic functions of $p$
variables. Jacobi's ``general theorem'' claimed that $x_0,x_1,\cdots,x_{p-1}$ were the roots of an algebraic
equation of degree $p$ whose coefficient were single-valued, $2p$-periodic functions of $u_0,u_1,\cdots,u_{p-1}$.
Therefore, the elementary symmetric functions of $x_0,x_1,\cdots,x_{p-1}$ could be expressed by means of
single-valued functions in $C^p$. In particular, Jacobi considered the case $p=2$ (\cite{Jacobi1881}, vol. 2,
7-16). His ideas were successfully developed by A. G\"{o}pel in 1847 (and, independently of him, J. G. Rosenhain
in 1851). The required 4-fold periodic functions of two complex variables were expressed as the ratio of two
$\theta$-series of two complex variables obtained by a direct and cumbersome computation.   This involved an
impressive amount of calculations and could hardly be extended to the case $p
> 2$. Following a completely different route Weierstrass was able
to solve the problem for any $p$. Because of his achievements he
was awarded a doctor degree {\em honoris causa} from the
K\"{o}nigsberg University, and two years later he was hired to
teach at the Berlin Gewerbeinstitut (later Gewerbeakademie, today
Technische Universit\"{a}t). Eventually, in the Fall of 1856
Weierstrass  was named {\em Extraordinarius} at the Berlin
University.

\section{Weierstrass' early papers} \label{section 1}
\setzero \vskip-5mm \hspace{5mm }

In the address he gave in 1857 upon entering the Berlin Academy,
Weierstrass recognized the ``powerful attraction'' which the
theory of elliptic functions had exerted on him since his student
days. In order to become a school teacher in 1839 Weierstrass had
entered the Theological and Philosophical Academy of M\"{u}nster,
where he attended for one semester Gudermann's lectures on
elliptic functions and became familiar with the concept of uniform
convergence which Gudermann had introduced in his papers in 1838.
Elliptic functions constituted the subject of Weierstrass' very
first paper, an essay he wrote  in autumn 1840 for obtaining his
{\em venia docendi}. His starting point was Abel's claim that the
elliptic function which is the inverse of the elliptic integral of
first kind (sn$u$ in the symbolism Weierstrass took from
Gudermann) could be expressed as the ratio of two convergent power
series of $u$, whose coefficients are entire functions of the
modulus of the integral. Weierstrass succeeded in proving that
sn$u$ (and similarly cn$u$  and dn$u$) could be represented as
quotient of certain functions, which he named {\em Al}-functions
in honor of Abel and which he was able to expand in convergent
power series.

Working in complete isolation without any knowledge of Cauchy's
related results, two years before Laurent, Weierstrass (1841)
succeeded in establishing the Laurent expansion of a function in
an annulus. In the paper he made an essential use of integrals,
and proved the Cauchy integral theorem for annuli (\cite{Weier94},
vol. 1, 51-66).  In a subsequent paper he stated and proved three
theorems on power series. Theorems A) and B) provided estimates
(Cauchy inequalities) for the coefficients of a Laurent series in
one (and several) complex variables, while Theorem C) was the
double series theorem nowadays called after him. As a consequence
of it Weierstrass obtained the theorem on uniform differentiation
of convergent series (\cite{Weier94}, vol. 1, 67-74). Apparently,
this paper marked a turning point in Weierstrass' analytic methods
for he gave up integrals and choose the power series approach to
treat the theory of function of one or more variables on a par.
This work was completed by a paper he wrote in spring 1842. There
Weierstrass proved that a system of $n$ differential equations
        $$\frac{dx}{dt}= G_i(x_1, \cdots, x_n) \mbox{    }(i = 1,\cdots, n)
  \mbox{   }(G_i(x_1, \cdots, x_n) \mbox{ polynomials})$$
can be solved by a system of $n$ unconditionally and uniformly
convergent power series satisfying  prescribed initial conditions
for $t$=0. In addition, he also showed how the power series
$$x_i=P_i(t-t_0, a_1, \cdots, a_n) \mbox{       }   (i = 1, \cdots, n;
t_0, a_1, \cdots, a_n \mbox{ fixed})$$ convergent in a disk centered at $t_0$ could be analytically continued
outside the disk. Thus, by the early 1840s the essential results of Weierstrass' approach to the theory of
analytic functions were already established. His papers, however, remained in manuscript and had no influence on
the contemporary development of mathematics.

\section{Abelian functions and integrals} \label{section 2}
\setzero \vskip-5mm \hspace{5mm }

Weierstrass' 1854 paper (\cite{Weier94}, vol. 1, 133-152) gave ``a short overview'' of the work on Abelian
functions which he had developed ``several years ago'' and summarized in the annual report of the Braunsberg
Gymnasium for 1848-49. Weierstrass began by considering the polynomial $R(x)= (x-a_0)(x-a_1)\cdots (x-a_{2n})$,
with $a_i$ real numbers satisfying the inequalities $a_i> a_{i+1}$. He decomposed $R(x)$ into the factors
$P(x)=\prod_{k=1}^n(x-a_{2k-1})$, $Q(x)=\prod_{k=1}^n(x-a_{2k})$ and considered the system
 \begin{eqnarray}\label{2.1}
   u_m=\sum_{j=1}^n \int_{a_{2j-1}}^{x_j}\frac{P(x)}{x-a_{2m-1}}\frac{dx}{2\sqrt{R(x)}} \mbox{     } \mbox{     }(m=1,\cdots,n).
  \end{eqnarray}
The task Weierstrass gave himself was to ``establish in detail'' Jacobi's theorem which he considered ``the
foundations of the whole theory''.  As Jacobi had remarked,  for given values of $x_1, x_2 ,\cdots, x_n$ the
quantities $u_1, u_2 ,\cdots ,  u_n$ have infinitely many different values. ``Conversely, if the values of $u_1,
u_2, \cdots, u_n$ are given, then the values of $x_1, x_2 , \cdots, x_n$ as well as the corresponding values of
$\sqrt{ R(x_1)},\sqrt{R(x_2)},\cdots,\sqrt{R(x_n)}$ are uniquely determined''. Moreover, ``$x_1, x_2 ,\cdots, x_n$
are roots of a (polynomial) equation of degree $n$ whose coefficient are completely determined, single-valued
functions of the variables $u_1, u_2 ,\cdots, u_n$''. Analogously, Weierstrass added, there exists a polynomial
function of $x$, whose coefficients are also single-valued functions of $u_1, u_2 ,\cdots,  u_n$ which gives the
corresponding values of $\sqrt{R(x_1)},\sqrt{R(x_2)},\cdots,\sqrt{R(x_n)}$  for $x = x_1, x_2 , \cdots,  x_n$.
Every rational symmetric function of $x_1, x_2 , \cdots, x_n$ could consequently be regarded as a single-valued
function of $u_1, u_2 , \cdots,  u_n$. Weierstrass considered the product $L(x) =(x-x_1)(x-x_2)\cdots(x-x_n)$ and
the $2n+1$ single-valued functions $Al(u_1,u_2,u_n)_m=\sqrt{h_m L(a_m)}$ $(m = 0,\cdots , 2n+1)$, with $h_m$
suitable constants, which he called Abelian functions, ``since they are the ones which completely correspond to
the elliptic functions'' to which they reduce when $n = 1$.

He was able to expand his $Al$-functions in convergent power series and, on the basis of Abel's theorem, he
succeeded in establish the ``principal property'' of such functions, i.e. an addition theorem according to which $
Al(u_1+v_1,u_2+v_2, \cdots,u_n+v_n)_m$ are rationally expressed in terms of $Al(u_1,u_2, \cdots,u_n)_m$,
$Al(v_1,v_2, \cdots,v_n)_m$ and their first-order partial derivatives. Eventually, he determined the algebraic
equation whose coefficients were expressed in terms of $Al$-functions, and whose roots were the quantities $x_1,
x_2 , \cdots, x_n$ satisfying equations \ref{2.1}) for arbitrary $u_1, u_2 ,\cdots,  u_n$. However, as Dirichlet
commented, in his paper Weiertrass ``gave only partial proofs of his results and lacked the intermediate
explanations'' (\cite{Dugac1973}, 52).

Two years later Weierstrass resumed this work and published in Crelle's \textit{Journal} the first part of an
expanded and detailed version of it (\cite{Weier94}, vol. 1, 297-355). As he had done in his 1854 paper,
Weierstrass considered the polynomial $R(x)=A(x-a_1)(x-a_2)\cdots(x-a_{2\rho+1})$, and the analogous product
$P(x)=\prod_{j=1}^{\rho}(x-a_j)$, $ (j=1,\cdots,{\rho})$ where this time the $a_j$ were any complex numbers such
that $a_j\neq a_k$ for $j\neq k$. Instead of equations \ref{2.1}) he considered the corresponding system of
differential equations
  \begin{equation}\label{2.2}
  du_m=\sum_{j=1}^\rho \frac{1}{2}\frac{P(x_j)}{x_j-a_m}\frac{dx}{\sqrt{R(x_j)}} \mbox{    }(m = 1,\cdots, \rho)
   \end{equation}
and formulated Jacobi inversion problem as the question to find solutions $x_j = x_j(u_1,\cdots, u_\rho)$ of the
system \ref{2.2}) satisfying the initial conditions $x_j(0,\cdots,0) = a_j$  $(j = 1,\cdots, \rho)$. In order to
obtain the elliptic functions as a special case ($\rho=1$), he also gave a slightly different form to his
$Al$-functions with respect to his previous paper.
    Weierstrass succeeded in proving that the solutions $x_j = x_j(u_1,\cdots, u_\rho)$ are
    single-valued functions
    of $u_1,\cdots, u_\rho$ in the neighborhood of the origin.
    They could be considered as the roots of a polynomial equation of degree $\rho$,
    whose coefficients were given in terms of $Al$-functions which, for any bounded value
    of ($u_1,\cdots, u_\rho$),
     are single-valued functions expressed as quotient of power series. Then, the symmetric
     functions of $x_j = x_j(u_1,\cdots, u_\rho)$
      have ``the character of rational functions''.
      From these results, however, Weierstrass was unable to show that each Abelian
      function could be represented as the ratio of two everywhere convergent power
      series. ``Here we encounter a problem that, as far as I know, has not yet been studied
      in its general form, but is nevertheless of particular importance for the theory of
      functions'' (\cite{Weier94}, vol. 1, 347).

      In the course of his life he returned many times to this problem
      in an attempt to solve it (see below). Even the factorization
      theorem
      for entire functions that Weierstrass was able to establish some 20 years later (see Section 5)
      can be regarded
      as an outcome of this research for it provided a positive answer to the problem in
      the case of one variable (\cite{Kolm96}, 247).
In order to show that his approach permitted one to treat the
theory of elliptic and Abelian function on a par, in the
concluding part of his 1856 paper Weierstrass presented a detour
on elliptic functions, where he summarized the main results he had
obtained in 1840. However, the promised continuation of
Weierstrass' paper never appeared. Instead, a completely new
approach to the theory of Abelian integrals was published by
Bernhard Riemann (1826-1866) in 1857 (\cite{Riem90}, 88-144) which
surpassed by far anything Weierstrass had been able to produce.

In the introductory paragraphs of his paper Riemann summarized the
geometric approach to complex function theory he had set up in his
1851 thesis (\cite{Riem90}, 3-48). There he defined a complex
variable $w$ as a function of $x+iy$ when $w$ varies according to
the equation $i\frac{\partial  w}{\partial x}= \frac{\partial
w}{\partial y}$ ``without assuming an expression of $w$ in terms
$x$ and $y$''. Accordingly, ``by a well known theorem'' - Riemann
observed without mentioning Cauchy - a function $w$ can be
expanded in a power series $\sum a_n(z-a)^n$ in a suitable disk
and ``can be continued analytically outside it in only one way''
(\cite{Riem90}, 88). For dealing with multi-valued functions such
as algebraic functions and their integrals Riemann introduced one
of his deepest achievements, the idea of representing the branches
of a function by a surface multiply covering the complex plane (or
the Riemann sphere).  Thus, ``the multi-valued function has only
{\em one} value defined at each point of such a surface
representing its branching, and can therefore be regarded as a
completely determined  (=single-valued) function of position on
this surface'' (\cite{Riem90}, 91). Having introduced such basic
topological concepts as crosscuts and order of connectivity of a
surface, Riemann could state the fundamental existence theorem of
a complex function on the surface, which he had proved in his
dissertation by means of a suitable generalization of the
Dirichlet principle. This theorem establishes the existence of a
complex function in terms of boundary conditions and behavior of
the function at the branch-points and singularities. Then Riemann
developed the theory of Abelian functions proper. It is worth
noting that, in spite of the fact that both Weierstrass and
Riemann gave their paper the same title and used the same wording,
they gave it a different meaning. Whereas the former defined
Abelian functions to be the single-valued, analytic functions of
several complex variables related to his solution of the Jacobi
inversion problem, the latter understood Abelian functions to be
the integrals of algebraic functions introduced by Abel's theorem.
In the first part of his paper Riemann developed a general theory
of such functions and integrals on a surface of any genus $p$,
"insofar as this  does not depend on the consideration of
$\theta$-series" (\cite{Riem90}, 100).  He was able to classify
Abelian functions (integrals) into three classes according to
their singularities, to determine the meromorphic functions on a
surface, and to formulate Abel's theorem in new terms, thus
throwing new light on the geometric theory of birational
transformations. The second part of the paper was devoted to the
study of $\theta$-series of $p$ complex variables, which express
``the Jacobi inverse functions of  $p$ variables for an arbitrary
system of finite integrals of equiramified, ($2p+1$)-connected
algebraic functions'' (\cite{Riem90}, 101). In this part Riemann
gave a complete solution of Jacobi inversion problem without
stating it as a special result. He regarded the work of
Weierstrass as a particular case, and mentioned the ``beautiful
results'' contained in the latter's 1856 paper, whose continuation
could show ``how much their results and their methods coincided''.
However, after the publication of Riemann's paper Weierstrass
decided to withdraw the continuation of his own. Even though
Riemann's work ``was based on foundations completely different
from mine, one can immediately recognize that his results coincide
completely with mine'' Weierstrass later stated
 (\cite{Weier94}, vol. 4, 9-10). ``The proof of this requires some
research of algebraic nature''. By the end of 1869 he had not been
able to overcome all the related ``algebraic difficulties''. Yet,
Weierstrass thought he had succeeded in finding the way to
represent any single-valued $2p$-periodic (meromorphic) function
as the ratio of two suitable $\theta$-series, thus solving the
general inversion problem. However, Weierstrass' paper
(\cite{Weier94}, vol. 2, 45-48) was flawed by some inaccuracies
that he himself later recognized in a letter to Borchardt in 1879
(\textit{ibid.}, 125-133). In particular, Weierstrass (mistakenly)
stated that any domain of ${C^{n}}$ is the natural domain of
existence of a meromorphic function. (This mistake was to be
pointed out in papers by F. Hartogs and E.E. Levi in the first
decade of the 20$^{th}$ century). By 1857, in his address to the
Berlin Academy Weierstrass limited himself to state that ``one of
the main problems of mathematics'' which he decided to investigate
was ``to give an actual representation'' of Abelian functions. He
recognized that he had published results ``in an incomplete
form''. ``However - Weierstrass continued - it would be foolish if
I were to try to think only about solving such a problem, without
being prepared by a deep study of the methods that I am to use and
without first practicing on the solution of less difficult
problems'' (\cite{Weier94}, vol. 1, 224). The realization of this
program became the scope of his University lectures.

\section{Weierstrass' lectures} \label{section 3}
\setzero \vskip-5mm \hspace{5mm }

In response to Riemann's achievements, Weierstrass devoted himself
to ``a deep study of the methods'' of the theory of analytic
functions which in his view provided the foundations of the whole
building of the theory of both elliptic and Abelian functions. As
Poincar\'{e} once stated, Weierstrass' work could be summarized as
follows: 1) To develop the general theory of functions, of one,
two and several variables. This was ``the basis on which the whole
pyramid should be built''. 2) To improve the theory of the
elliptic functions and to put them into a form which could be
easily generalized to their ``natural extension'',  the Abelian
functions. 3) Eventually, to tackle the Abelian functions
themselves.

Over the years the aim of establishing the foundations of analytic
function theory with absolute rigor on an arithmetic basis became
one of Weierstrass' major concerns.  From the mid-1860s to the end
of his teaching career Weierstrass used to present the whole of
analysis in a two-year lecture cycle as follows:
\begin{enumerate}
\item  Introduction to analytic function theory,

\item  Elliptic functions,

\item  Abelian functions,

\item Applications of elliptic functions or, alternatively,
Calculus of variations.
\end{enumerate}

All of these lectures, except for the introduction to analytic
function theory, have been published in Weierstrass' {\em Werke}.
For some twenty years he worked out his theory of analytic
functions through continuous refinements and improvements, without
deciding to publish it himself. Weierstrass used to present his
discoveries in his lectures, and only occasionally communicated
them to the Berlin Academy. This attitude, combined with his
dislike of publishing his results in printed papers and the fact
that he discouraged  his students from publishing lecture notes of
his courses, eventually gave Weierstrass' lectures an aura of
uniqueness and exceptionality.

\section{Conversations in Berlin} \label{section 4}
\setzero \vskip-5mm \hspace{5mm }

In the Fall of 1864, when Riemann was staying in Pisa because of
his poor health conditions, the Italian mathematician F. Casorati
travelled to Berlin to meet Weierstrass and his colleagues. Rumors
about new discoveries made by Weierstrass, combined with lack of
publications, motivated Casorati's journey.

 ``Riemann's things are creating difficulties in Berlin'', Casorati
 recorded in his notes.  Kronecker claimed that ``mathematicians $\cdots$ are a bit
 arrogant ({\em hochm\"{u}tig}) in using the concept of
 function''. Referring to Riemann's proof of the Dirichlet
 principle, Kronecker remarked that Riemann himself, ``who is generally very precise,
 is not beyond censure in this regard'' (\cite{Bottazzini1986}, 262).

Kronecker added that in Riemann's paper on Abelian functions the
$\theta$-series in several variables ``came out of the blue''.
Weierstrass claimed that ``he understood Riemann, because he
already possessed the results of his [Riemann's] research''. As
for Riemann surfaces, they were nothing other than ``geometric
fantasies''. According to Weierstrass, ``Riemann's disciples are
making the mistake of attributing everything to their master,
while many [discoveries] had already been made by and are due to
Cauchy, etc.; Riemann did nothing more than to dress them in his
manner for his convenience''. Analytic continuation was a case in
point. Riemann had referred to it in various places but, in
Weierstrass's and Kronecker's opinion, nowhere he had treated it
with the necessary rigor. Weierstrass observed that Riemann
apparently shared the idea that it is always possible to continue
a function to any point of the complex plane along a path that
avoids critical points (branch-points, and singularities). ``But
this is not possible'', Weierstrass added. ``It was precisely
while searching for the proof of the general possibility that he
realized it was in general impossible''. Kronecker provided
Casorati with the example of the (lacunary) series

  \begin{equation}\label{4.1}
 \theta_0(q) = 1 + 2 \sum_{n\geq 1}q^{n^2}
  \end{equation}
which is convergent for $|q | < 1$, and has the unit circle as a
natural boundary. Its unit circle is ``entirely made of points
where the function is not defined, it can take any value there'',
Weierstrass observed. He had believed that points in which a
function ``ceases to be definite'' - as was the case of the
function $e^{1/x}$ at $x = 0$ because ``it can have any possible
value'' there - ``could not form a continuum, and consequently
that there is at least one point $P$ where one can always pass
from one closed portion of the plane to any other point of it''.
$\theta_0(q)$ provided an excellent example of this unexpected
behavior. This series also played a significant role in
Weierstrass' counter-example of a continuous nowhere
differentiable function (see Section 6).

\section{Further criticism of Riemann's methods} \label{section 5}
\setzero \vskip-5mm \hspace{5mm }

Apparently, Riemann's theory of complex functions seems to have
been the background of Weierstrass' work and lectures. Evidence of
this is provided by his (unpublished) correspondence with his
former student H. A. Schwarz from 1867 up to 1893. One of the
first topics they discussed was Riemann mapping theorem. In his
thesis Riemann had claimed that ``two given simply connected plane
surfaces can always be mapped onto one another in such a way that
each point of the one corresponds to a unique point of the other
in a continuous way and the correspondence is conformal; moreover,
the correspondence between an arbitrary interior point of the one
and the other may be given arbitrarily, but when this is done the
correspondence is determined completely'' (\cite{Riem90}, 40).
Riemann's proof of the mapping theorem rested on a suitable
application of  Dirichlet's principle. Because of his criticism of
this principle, in Weierstrass' view the Riemann mapping theorem
remained a still-open question, worthy of a rigorous answer.
Following Weierstrass' suggestion, Schwarz tackled this question
after his student days and succeeded in establishing the theorem
in particular cases, without resorting to the questionable
principle. In a number of papers he gave the solution of the
problem of the conformal mapping of an ellipse - or, more
generally, of a plane, simply connected figure, with boundaries
given by pieces of analytic curves which meet to form non-zero
angles - onto the unit disk, by using suitable devices as the
lemma and the reflection principle, both named after him
(\cite{Schw90}, vol. 2, 65-132).

In 1870 Schwarz discovered his alternating method. ``With this
method - he stated by presenting it in a lecture - all the
theorems which Riemann has tried to prove in his papers by means
of the Dirichlet principle, can be proved rigorously''
(\cite{Schw90}, vol. 2, 133). He submitted to Weierstrass an
extended version of the paper, and in a letter of July 11, 1870
Schwarz asked him whether he had ``objections to raise''.
Apparently, Weierstrass' answer has been lost. It is quite
significant, however, that three days later, on July 14, 1870
Weierstrass presented to the Berlin Academy his celebrated
counterexample to the Dirichlet principle (\cite{Weier94}, vol. 2,
49-54), and then submitted Schwarz's 1870 paper for publication in
the {\em Monatshefte} of the Academy. Two years later, in a letter
of June 20, 1872 Schwarz called Weierstrass' attention to the
still widespread idea that a continuous function always is
differentiable. As the French mathematician Joseph Bertrand had
made this claim in the opening pages of his {\em Trait\'{e}},
Schwarz ironically wondered about asking Bertrand to prove that
 \begin{equation}\label{5.1}
 f(x) = \sum_{n\geq 1}\frac{\sin {n^{2}x}}{{n^2}}
  \end{equation}
has a derivative. One month later, on July 18, 1872 Weierstrass presented the Academy with his celebrated example
of a continuous, nowhere differentiable function
 \begin{equation}\label{5.2}
 f(x) = \sum_{n\geq 0}{b^{n}\cos{a^{n}}x\pi}
  \end{equation}
where $a$= is an odd integer, $0<b<1$, and $ab>1+ \frac{3}{2}\pi$. According to Riemann's students, Weierstrass
remarked, the very same function  \ref{5.1}) mentioned by Schwarz had been presented by Riemann in 1861 or perhaps
even earlier in his lectures as an example of continuous nowhere differentiable function. ``Unfortunately
Riemann's proof has not been published'', Weierstrass added, and ``it is somewhat difficult to prove'' that 5.1)
has this property, he concluded before producing his own example (\cite{Weier94}, vol. 2, 71-74).

Only by the end of 1874 was Weierstrass able to overcome a major
difficulty which for a long time had prevented him from building a
satisfactory theory of single-valued functions of one variable.
This was the proof of the representation theorem of a
single-valued function as a quotient of two convergent power
series. As he wrote on the same day (December 16, 1874) to both
Schwarz and S. Kovalevskaya, this was related to the following
question: given an infinite sequence of constants $\{a_n\}$ with
$\lim {|a_n|} = \infty$ does there always exist an entire,
transcendental function $G(x)$ which vanishes at $\{a_n\}$ and
only there? He had been able to find a positive answer to it by
expressing $G(x)$ as the product$\prod_{n\geq 1}{E(x,n)}$ of
"prime functions"

$E(x,0)=1+x$,

$\cdots$

$E(x,n)$=$(1+x)\exp({\frac{x}{1}}+{\frac{x^{2}}{2}}+\cdots+{\frac{x^{n}}{n})}$

\noindent which he introduced there for the first time. The ``until now only conjectured'' representation theorem
followed easily. This theorem constituted the core of  Weierstrass's 1876 paper on the ``systematic foundations''
of the theory of analytic functions of one variable (\cite{Weier94}, vol. 2, 77-124). In spite of his efforts,
however, he was not able to extend his representation theorem to single-valued functions of several variables.
``This is regarded as unproved in my theory of Abelian functions'' Weierstrass admitted in his letter to
Kovalevskaya. (For 2 variables this was done by Poincar\'{e} in 1883 and later extended by Cousin in 1895
following different methods from Weierstrass'). Four days later Weierstrass wrote to Schwarz stating that
Riemann's (and Dirichlet's) proof of Cauchy integral theorem by means of a double integration process was in his
opinion not a "completely methodical" one. On the contrary, a rigorous proof could be obtained by assuming the
fundamental concept of analytic element (and its analytic continuation) and by resorting to Poisson integral for
the disk, as Schwarz himself had shown in his paper on the integration of the Laplace equation. Criticism of
Riemann's ideas and methods were also occasionally expressed by Weierstrass
    in his letters to Kovalevskaya \cite{Mitt23}. On August 20, 1873 he
    was pleased to quote an excerpt from a letter of Richelot to himself ``in which a decisive
    preference was expressed for the route chosen by Weierstrass in the theory of Abelian
    functions as opposed to Riemann's and Clebsch's''. On January 12, 1875 Weierstrass announced
    to Kovalevskaya his intention of presenting the essentials of his approach to Abelian
    functions in a series of letters to Richelot where he hoped ``to point out the uniqueness
    of my method without hesitation and to get into a criticism of Riemann and Clebsch''.

Weierstrass openly stated his criticism of Riemann's methods in a
often-quoted ``confession
    of faith'' he produced to Schwarz on October 3, 1875: ``The more I think about the principles
    of function theory - and I do it incessantly - the more I am convinced that this must be built
    on the basis of algebraic truths, and that it is consequently not correct when the `transcendental',
    to express myself briefly, is taken as the basis of simple and fundamental algebraic
    propositions. This view seems so attractive at first sight, in that through it Riemann
    was able to discover so many of the important properties of  algebraic functions''. Of course,
    Weierstrass continued, it was not a matter of methods of discovery. It was ``only a matter of systematic foundations'' (\cite{Weier94}, vol. 2, 235).
    It is worth remarking that Weierstrass added he had been ``especially strengthened [in his
    belief] by his continuing study of the theory of analytic functions of several variable''.

\section{Weierstrass' last papers} \label{section 6}
\setzero \vskip-5mm \hspace{5mm }

After Mittag-Leffler, Poincar\'{e} and Picard had deeply extended
the results of his 1876 paper following ``another way'' different
from his own, Weierstrass felt it necessary to explain his
approach to complex function theory and to compare it with those
of Cauchy and Riemann. He did this in a lecture that he delivered
at the Berlin Mathematical Seminar on May 28, 1884 \cite{Weier23}
. Even though ``much can be done more easily by means of Cauchy's
theorem'', Weierstrass admitted, he strongly maintained that the
general concept of a single-valued analytic function had to be
based on simple, arithmetical operations. His discovery of both
continuous nowhere differentiable functions and series having
natural boundaries strengthened him in this view. ``All
difficulties vanish'', he stated, ``when one takes an arbitrary
power series as the foundation of an analytic function''
(\cite{Weier23}, 3).

Having summarized the main features of his own theory, including
in particular the method of analytic continuation,  he advanced
his criticism of Riemann's general definition of a complex
function (\textit{see} Section 2). This was based on the existence
of first-order partial derivatives of  functions of two real
variables, whereas ``in the current state of knowledge'' the class
of functions having this property could not be precisely
delimited. Moreover, the existence of partial derivatives required
an increasing number of assumptions when passing from one to
several complex variables. On the contrary, Weierstrass concluded,
his own theory could ``easily'' be extended to functions of
several variables.

A major flaw in Riemann's concept of a complex function had been
discovered and published by Weierstrass in 1880. The main theorem
of his paper stated  that a series of rational functions,
converging uniformly inside a disconnected domain may represent
different analytic functions on disjoint regions of the domain
(\cite{Weier94}, vol. 2, 221). Thus, Weierstrass commented, ``the
concept of a monogenic function of a complex variable does not
coincide completely with the concept of dependence expressed by
(arithmetic) operations on quantities'', and in a footnote he
pointed out that ``the contrary statement had been made by
Riemann'' in his thesis. Before proving his theorem Weierstrass
discussed an example he had expounded in his lectures ``for many
years''. By combining the theory of linear transformations of
elliptic $\theta$-functions with the properties of the lacunary
series \ref{4.1}), Weierstrass was able to prove that the series
 \begin{equation}
 F(x) = \sum_{n\geq 0}\frac {1} {x^{n}+x^{-n}}
  \end{equation}
is convergent for $|x | < 1$, and $|x | > 1$, but  ``in each region of its domain of convergence it represents a
function which cannot be continued outside the boundary of the region'' (\cite{Weier94}, vol. 2, 211). (It is
worth noting that $1+4F(x)$ = $\theta^{2}_0(x) $).

This remark allowed Weierstrass to clarify an essential point of
function theory, which deeply related the problem of the
analytical continuation of a complex function to the existence of
real, continuous nowhere differentiable functions. In order to
explain this relation Weierstrass considered the series
$\sum_{n\geq 0}b^{n}x^{a^n}$ which is absolutely and uniformly
convergent in the compact disk $|x | \leq 1$, when $a$ is an odd
integer, $0 < b < 1$. By a suitable use of his example of a
continuous nowhere differentiable function \ref{5.2}), he
concluded that under the additional condition $ab>1+
\frac{3}{2}\pi$  the circle $|x | = 1$ reveals to be the natural
boundary of the series. Contrary to his habit, in 1886 Weierstrass
reprinted this paper in a volume which collected some of his last
articles, including a seminal paper where he stated his celebrated
``preparation theorem'' together with other theorems on
single-valued functions of several variables that he used to
expound in his lectures on Abelian functions \cite{Weier86}.

\section{Conclusion}
\setzero \vskip-5mm \hspace{5mm }

From the 1840s to the end of his life Weierstrass continued to
study the theory of Abelian functions, devoting an incredible
amount of work to the topic. This theory was the background of
many of the results he presented in his papers and lectures, or
discussed in his letters to colleagues. In spite of his efforts,
however, Weierstrass never succeeded in giving it the complete,
rigorous treatment he was looking for. The huge fourth volume of
his \textit{Mathematishe Werke} (published posthumously) collects
the lectures on Abelian functions he gave in Winter semester
1875-76 and Summer semester 1876. Two thirds of it is devoted to
algebraic functions and Abelian integrals, and only the remaining
one third to the (general) Jacobi inversion problem. Thus, the
editors of the volume, Weierstrass' former students G. Hettner and
J. Knoblauch, could aptly state in the preface that the theory of
Abelian functions (in Weierstrass' sense) ``is sketched only
briefly'' there. It was not an irony of the history if Weierstrass
failed in his pursuit of his his main mathematical goal whereas
the machinery that he created to attain it in response to
Riemann's ``geometric fantasies'' became an essential ingredient
of modern analysis. The contradiction between Weierstrass'
approach, in which all geometric insight was lacking, and
Riemann's geometric one remained effective until the early decades
of  the 20$^{th}$ century, when the theory of functions of several
complex variables began to be established in modern terms.

\end{document}